\documentclass{article}
\usepackage{amssymb, verbatim}
\usepackage{hyperref}
\usepackage{amsmath}
\usepackage{tikz-cd}

\def\Ric{\mathop{\rm Ric}}
\def\cRic{\mathop{\rm R\i\makebox[0pt]{\raisebox{5pt}{\tiny$\circ$\;\,}}c}}

\def\dist{\mathop{\rm dist}}

\def\Riem{\mathop{\rm Rm}}

\def\Vol{\mathop{\rm Vol}}

\def\RRR{\mathop{\mathcal{R}}}

\def\CP{\mathop{\mathbb{CP}}}

\def\be{\begin{eqnarray}}
\def\ee{\end{eqnarray}}
\def\beg{\begin{eqnarray*}}
\def\ees{\end{eqnarray*}}

\def\XXint#1#2#3{{\setbox0=\hbox{$#1{#2#3}{\int}$ }
\vcenter{\hbox{$#2#3$ }}\kern-.5\wd0}}

\newcommand{\qed}{\hfill$\Box$}
\newtheorem{theorem}{Theorem}[section]

\newtheorem{lemma}[theorem]{Lemma}

\newtheorem{corollary}[theorem]{Corollary}

\date{\small\it September 9, 2013}
\author{Brian Weber}
\title{Harnack Inequalities for Critical 4-Manifolds with a Ricci Curvature Bound}

\begin{document}

\maketitle

\section{Introduction}

We develop certain elliptic inequalities for metrics on 4-manifolds with a lower bound on Ricci curvature and some form of elliptic system; specifically we prove two Harnack-type theorems on low-energy regions of such manifolds, and explore some consequences.
Our results are of collapsing type: volume ratios and Sobolev constants are immaterial.
First consider the $s$-local curvature radius
\be
r^s_{\mathcal{R}}(p)\;=\;\sup\,\left\{\;r\,\in\,(0,\,s)\;\big|\;|\Riem|<r^{-2}\;{\rm on}\;B(p,\,r)\;\right\},
\ee
(also called the local curvature radius with cutoff $s$); we use $r_{\mathcal{R}}$ for $r_{\mathcal{R}}^\infty$.
Besides being a local bound on sectional curvature, the basic role of this function on critical Riemannian manifolds has been made clear in \cite{CT} and elsewhere, where it was used in building N-structures on low energy regions of a Riemannian manifolds (see Section \ref{SectionPriorResults}).

The thrust of this paper is to obtain elliptic-type estimates governing the behavior of $r_{\mathcal{R}}$.
Unlike the typical proofs from the $L^p$ theory which rely on (here unavailable) Sobolev inequalities, our estimates hinge on geometric and topological considerations.
The elliptic system itself plays a background role; it is needed to enforce low-energy collapsing, as described in \cite{CT}, and then to ensure local convergence in at least the $C^2$ sense---continuous convergence of Riemann tensors.
The central argument, contained in the proof of Theorem \ref{ThmHarnackI}, is a contradiction involving the ``geometric vanishing'' theorem of \cite{Web2} and a blow-up argument on collapsing regions of Riemannian manifolds, where $C^2$ convergence implies a non-flat limit.

We have the following assumptions throughout.
First, our manifolds $(N^4,g)$ have critical metrics (see for example \cite{And4}), and in particular have an elliptic system.
Common examples of such metrics are Einstein, half conformally-flat with constant scalar curvature, Bach-flat with constant scalar curvature, and extremal K\"ahler metrics.
In all cases we assume a uniform lower bound on Ricci curvature: $\Ric\ge-\Lambda^2$.
In the constant scalar curvature (CSC) Bach-flat case (which includes the half-conformally flat case) and in the CSC K\"ahler case, this is our only assumption.
In the extremal K\"ahler case, we must assume that the gradient of the scalar curvature is bounded: if $X=\nabla{R}$ where $R$ is scalar curvature, then $|X|<\Lambda^3$.
In the case of any other other elliptic system, we actually require two-sided bounds on Ricci curvature: $|\Ric|\le\Lambda^2$.
The symbol $\Omega$ will indicate a pre-compact domain in $N^4$.

Our primary result is a Harnack-type estimate on $r_{\mathcal{R}}$; our other results are corollaries.
Of course the function $r_{\mathcal{R}}$ is controlled from below on the {\it interior} of $B(p,r_{\mathcal{R}}(p))$, but in general has no lower bounds on the ball's boundary or exterior.
The general inequality is
\be
r_{\mathcal{R}}(p')\;\ge\;r_{\mathcal{R}}(p)-\dist(p,p').
\ee
which implies $r_{\mathcal{R}}$ has Lipschitz constant 1.
The first theorem provides, in regions of small energy, definite lower control of $r_{\mathcal{R}}$ up to and beyond the boundary of $B(p,r_{\mathcal{R}}(p))$: in the interiors of low-energy domains, the largest values of $\left(r_{\mathcal{R}}\right)^{-1}$ are bounded in terms of its nearby lowest values.

\begin{theorem}[Harnack inequality I] \label{ThmHarnackI}
Given $K<\infty$, there exists $\epsilon_0=\epsilon_0(K)>0$ and $\delta_0=\delta_0(K)>0$ so that if $r_{\mathcal{R}}(p)<K\Lambda^{-1}$ and $\int_{B(p,\,2Kr_{\mathcal{R}}(p))}|\Riem|^2\;\le\;\epsilon_0$, then $r_{\mathcal{R}}\;\ge\;\delta_0\,r_{\mathcal{R}}(p)$
on $B(p,Kr_{\mathcal{R}}(p))$.
\end{theorem}
\begin{theorem}[Harnack inequality II] \label{ThmHarnackII}
Given $K<\infty$, $k>0$, and $\mu\in(0,1]$, there exist numbers $\epsilon_0=\epsilon_0(K)>0$ and $C=C(\mu,k,K)<\infty$ with the following.
If $r_{\mathcal{R}}(p)<K\Lambda^{-1}$ and  $\int_{B(p,\,2{r}_{\mathcal{R}}(p))}|\Riem|^2\;\le\;\epsilon_0$, then
\be
r_{\mathcal{R}}(p)^{-2k}\;\le\;\frac{C}{\Vol\,B(p,\,\mu{r}_{\mathcal{R}}(p))}\int_{B(p,\,\mu{r}_{\mathcal{R}}(p))}|\Riem|^{k}. \label{IneqHarnIIa}
\ee
\end{theorem}
To explain the relevance of the second theorem, note that for an arbitrary metric, $|\Riem|\le{r}_{\mathcal{R}}(p)^{-2}$ on $B(p,r_{\mathcal{R}}(p))$.
Theorem \ref{ThmHarnackII} says the reverse inequality holds in an average sense, provided energy is small.
Simply using this estimate on $r_{\mathcal{R}}^{-2}$ first, we also obtain
\be
r_{\mathcal{R}}(p)^{-2k}\;\le\;C\left(\frac{1}{\Vol\,B(p,\,\mu{r}_{\mathcal{R}}(p))}\int_{B(p,\,\mu{r}_{\mathcal{R}}(p))}|\Riem|^2\right)^{\frac{k}{2}}. \label{IneqHarnIIb}
\ee
These ``Harnack'' inequalities will likely be of further analytic interest.
The main consequence we prove is the following elliptic-type estimate on the curvature scale.
\begin{corollary}[Local elliptic estimates for the curvature radius] \label{CorEllipticEsts}
Given $K<\infty$ and $l\in\mathbb{N}$, there is an $\epsilon_0=\epsilon_0(K)>0$ and $C=C(K,l)$ so that if $r_{\mathcal{R}}(p)<K\Lambda^{-1}$ and
\be
\int_{B(p,\,2r_{\mathcal{R}}(p))}|\Riem|^2\;\le\;\epsilon_0,
\ee
then
\be
\sup_{B(p,r_{\mathcal{R}}(p))}\left|\nabla^lr_{\mathcal{R}}\right|\;\le\;C\,\left(r_{\mathcal{R}}(p)\right)^{1-l}.
\ee
\end{corollary}
If $l=1$ we see that $|\nabla{r}_{\mathcal{R}}|$ has an absolute bound---this is expected; we already noted $r_{\mathcal{R}}$ is Lipschitz.
When $l>1$, $\mu\in(0,1]$ and $\int_{B(p,2r_{\mathcal{R}}(p))}|\Riem|^2<\epsilon_0$, then in combination with (\ref{IneqHarnIIa}) and (\ref{IneqHarnIIb}) we have
\be
\begin{aligned}
\sup_{B(p,r_{\mathcal{R}}(p))}\left|\nabla^lr_{\mathcal{R}}\right|&\;\le\;\frac{C}{\Vol\,B(p,\,\mu{r}_{\mathcal{R}}(p))}\int_{B(p,\mu{r}_{\mathcal{R}}(p))}|\Riem|^{\frac{l-1}{2}} \\
\sup_{B(p,r_{\mathcal{R}}(p))}\left|\nabla^lr_{\mathcal{R}}\right|&\;\le\;C\left(\frac{1}{\Vol\,B(p,\,\mu{r}_{\mathcal{R}}(p))}\int_{B(p,\mu{r}_{\mathcal{R}}(p))}|\Riem|^2\right)^{\frac{l-1}{4}}
\end{aligned}
\ee
or the familiar-looking
\be
\begin{aligned}
\sup_{B(p,r_{\mathcal{R}}(p))}\left|\nabla^lr_{\mathcal{R}}\right|
&\;\le\;C\,\left(r_{\mathcal{R}}(p)\right)^{3-l}\left(\frac{1}{\Vol\,B(p,\,\mu{r}_{\mathcal{R}}(p))}\int_{B(q,\mu{r}_{\mathcal{R}}(q))}|\Riem|^2\right)^{1/2}.
\end{aligned}
\ee

{\bf Remark}.
In Theorem \ref{ThmHarnackI}, it would be nice to have an explicit approximation for $\delta_0$; possibly exponential decrease with increasing $K$.
As energy $\epsilon_0$ approaches zero, it is likely that $\delta_0$ approaches unity

{\bf Remark}.
As in the K\"ahler case, the half-conformally flat case does not actually require constant scalar curvature, but just a bound on the derivatuve of curvature: $|\nabla{R}|<\Lambda^3$.
But this is not a particularly natural assumption.

{\bf Remark}. 
In the extremal K\"ahler and CSC half-conformally flat cases, we conjecture that all conclusions hold without the $\Ric\ge-\Lambda^2$ assumption, provided uniform bounds on scalar curvature exist: $|R|<\Lambda^2$.
For other metrics with an elliptic system, such as CSC Bach-flat or harmonic curvature metrics (see \cite{Derd}), the situation is less clear.

\section{Prior Results} \label{SectionPriorResults}

In this section we outline the results and definitions that will be important for us.
The first is ``standard'' $\epsilon$-regularity; the second is a criterion for collapsing; the third is a criterion for flatness for certain 4-manifolds.
The last is the (rather involved) theory of N-structures.

\subsection{Epsilon-regularity, collapsing, and Ricci-pinched manifolds}

\begin{lemma}[Standard epsilon regularity] \label{LemmaStdEpsReg}
Given $K>0$ there exists numbers $\epsilon_0=\epsilon_0(K)>0$, $C=C(K)<\infty$ so that $r\le{K}\Lambda^{-1}$ and
\be
\frac{1}{Vol\,B(p,r)}\int_{B(p,r)}|\Riem|^2\;\le\;\epsilon_0\,r^{-4} \label{IneqEpsRegAssumption}
\ee
imply
\be
\sup_{B(p,r/2)}|\Riem|\;\le\;C\left(\frac{1}{Vol\,B(p,r)}\int_{B(p,r)}|\Riem|^2\right)^{\frac12}.
\ee
\end{lemma}
Lemma \ref{LemmaStdEpsReg} is normally used in a non-collapsing setting, for obvious reasons: if one assumes bounded volume ratios, say $\Vol\,B(p,r)\ge\delta{r}^n$, then one may measure $\int_{B(p,r)}|\Riem|^2$ against the a priori controlled quantity $r^{-4}\Vol\,B(p,r)$.
An argument found in \cite{CT}, effectively a contrapositive, extends its usefulness to the collapsing setting, essentially providing a way of forcing {\it collapse with locally bounded curvature}.
The significance of collapse with bounded (or locally bounded) curvature is explained in Section \ref{SubSecCollapsingStructs} below.
\begin{lemma}[Low energy collapse] \label{LemmaChTiForcedCollpasing}
Given $K>0$, $\tau>0$, there is an $\epsilon=\epsilon(\tau,K)>0$ so that $r_{\mathcal{R}}(p)\le{K}\Lambda^{-1}$ and $\int_{B(p,2r_{\mathcal{R}}(p))}|\Riem|^2\;\le\;\epsilon$ imply $\Vol\,B(p,\,r_{\mathcal{R}}(p))\;\le\;\tau\cdot{r}_{\mathcal{R}}(p)^{4}$.
\end{lemma}
\underline{\sl Pf}.
There is a point $q\in{B}(p,r_{\mathcal{R}}(p))$ with $|\Riem(q)|=r_{\mathcal{R}}(p)^{-2}$.
Now, assuming that $r_{\mathcal{R}}(p)^{-4}\,Vol\,B(p,r_{\mathcal{R}}(p))>\tau$, then choosing $\epsilon_0$ small enough we have (\ref{IneqEpsRegAssumption}).
But then the conclusion of Lemma \ref{LemmaStdEpsReg} holds, so
\be
r_{\mathcal{R}}(p)^{-2}\;\le\;C\left(\frac{1}{Vol\,B(p,2r_{\mathcal{R}}(p))}\int_{B(p,2r_{\mathcal{R}}(p))}|\Riem|^2\right)^\frac12.
\ee
Thus $r_{\mathcal{R}}(p)^{-4}Vol\,B(p,r_{\mathcal{R}}(p))<C^2\epsilon_0$.
Possibly choosing $\epsilon_0$ still smaller, we again have $r_{\mathcal{R}}(p)^{-4}Vol\,B(p,r_{\mathcal{R}}(p))\le\tau$.
\qed

The third result, from \cite{Web2}, is a flatness (``geometric vanishing'') theorem for certain manifolds with a Killing field.
We shall use it to conclude that certain blow-ups of collapsing manifolds, which are a priori non-flat, are in fact flat.
Specifically, if $\Omega$ is any domain in an $n$-manifold $N^n$, set
\be
M_X^{\Omega,s} \;=\;
\frac{{\sup}_{p\in{N}\setminus\Omega}\,\left\{\,|X(q)|\;\big|\;q\in{B}(p,\,r^s_{\mathcal{R}}(p))\;\setminus\;\Omega\right\}}{{\inf}_{p\in{N}\setminus\Omega}\,\left\{\,|X(q)|\;\big|\;q\in{B}(p,\,r^s_{\mathcal{R}}(p))\;\setminus\;\Omega\right\}}.
\ee
If $M_X^{\Omega,s}(p)<\infty$, we say $X$ has bounded $s$-local variation outside $\Omega$.
Then set
\be
M_X^{\infty}(p)\;=\;\inf_{s>0}\,M_X^{B(p,s),s}
\ee
If $M_X^{\infty}(p)<\infty$, we say $X$ has {\it asymptotically bounded local variation}.
Clearly $M_X^{\infty}(p)$ is independent of $p$, and we write simply $M_X^\infty$.
\begin{lemma}[Geometric Vanishing] \label{LemmaFlatMfld}
Assume $(N^4,\,g)$ is a complete 4-manifold with $\Ric\ge0$ and a nowhere-zero Killing field $X$ with $M_X^\infty<\infty$.
Then $N$ is flat provided:
\begin{itemize}
\item[{\it{a}})] The trace-free Ricci tensor is sufficiently pinched: $|\cRic|^2\le\frac{1}{24}R^2$; for instance if $(N^4,g)$ is Einstein
\item[{\it{b}})] The metric has zero scalar curvature and is half-conformally flat
\item[{\it{c}})] The metric is K\"ahler with respect to some complex structure, and has zero scalar curvature
\end{itemize}
\end{lemma}
The usefulness of this theorem is that, in the case of collapse with bounded curvature, we have N-structures whose associated locally-defined Killing fields have bounded local variation.
By passing to appropriate covers, we obtain complete manifolds where the Killing field(s) obtained from the N-structure automatically have asymptotically bounded local variation.

\subsection{Collapsing with bounded curvature: F- and N-structures} \label{SubSecCollapsingStructs}

The F- and N-structures of Cheeger-Gromov \cite{CG1} \cite{CG2} and Cheeger-Gromov-Fukaya \cite{CFG} will be decisive, so we define them precisely.
A number of variant definitions are available; ours is similar to the definition found in \cite{CR}, with the one main difference explained below.

An N-structure $\mathfrak{N}$ is a triple $(\Omega,\,\mathcal{N},\,\iota)$ where $\Omega$ is a domain in a differentiable manifold, $\mathcal{N}$ is a sheaf of nilpotent Lie algebras on $\Omega$, and $\iota:\mathcal{N}\rightarrow\mathcal{X}(\Omega)$ (called the action) is a sheaf monomorphism from $\mathcal{N}$ into the Lie algebra sheaf $\mathcal{X}(\Omega)$ of differentiable vector fields on $\Omega$, so that a collection of sub-structures $\mathcal{A}=\{(\mathcal{N}_i,\Omega_i,\iota_i)\}_i$ exists that satisfies the three conditions below.
In what follows, if $p\in\Omega_i$ its $\mathcal{N}_i$-stalk will be denoted $\mathcal{N}_{i,p}$ and its $\mathcal{N}$-stalk will be denoted $\mathcal{N}_p$.
\begin{itemize}
\item[{\it{i}})] (Completeness of the cover) The collection $\{\Omega_i\}$ is a locally finite cover of $\Omega$, and given $p\in\Omega$ there is at least one $\Omega_i$ so that $\mathcal{N}_{i,p}=\mathcal{N}_p$.
\item[{\it{ii}})] (Uniformity of the action) The lifted sheaf $\widetilde{\mathcal{N}_i}$ over the universal cover $\widetilde{\Omega_i}\rightarrow{\Omega_i}$ is a constant sheaf (each stalk is canonically isomorphic to the Lie algebra of global sections $\widetilde{\mathcal{N}_i}(\widetilde{\Omega_i})$).

As an aside, the lifted action $\widetilde{\iota_i}:\widetilde{\mathcal{N}_i}\rightarrow\mathcal{X}(\widetilde{\Omega_i})$ is not uniquely defined but depends on a choice of fundamental domain.
This manifests on the $\Omega_i$ as a holonomy phenomenon on stalks.
\item[{\it{iii}})] (Integrability of the action) Given $\Omega_i$, there is a connected, simply-connected nilpotent Lie group $G_i$ so that for any choice of $\widetilde{\iota_i}$ there is an action of $G_i$ on $\widetilde{\Omega_i}$ whose derived action is equal to the image of the Lie algebra of sections $\widetilde{\mathcal{N}_i}(\widetilde{\Omega_i})$ under $\widetilde{\iota_i}$.
\end{itemize}
An N-structure is an called an F-structure the associated sheaf $\mathcal{N}$ is abelian.
The difference between our definition of F-structures and the common definition is that we do not require that a torus acts on a finite normal cover of $\Omega_i$, but that some $\mathbb{R}^k$ acts on its universal cover.
This is a convenience in that we will make frequent passages to universal covers, and wish to refer to the structures obtained there as F- or N-strucutres, whether orbits are bounded or not.

Via the action of the groups $G_i$ on covers, an N-structure partitions $\Omega$ into orbits; the orbit through $p\in\Omega$ is denoted $\mathcal{O}_p$.
An orbit $\mathcal{O}_p$ is called {\it singular} if its dimension is not equal to the dimension of the stalk $\mathcal{N}_p$.
In addition, orbits may be {\it exceptional}; these are orbits for which nearby orbits are identified to it in a $k$-to-$1$ fashion.
An example would be $\mathbb{S}^3\subset\mathbb{C}^2$ with a Killing field given by differentiating the action $t\mapsto(e^{2\pi{i}t/k}z_1,\,e^{2\pi{i}lt/k}z_2)$, $k,l\in\mathbb{Z}$ relatively prime, so the two exceptional orbits are the points of the form $(z_1,0)$ and $(0,z_2)$ in $\mathbb{S}^3$.

The {\it rank} of an N-structure $\mathfrak{N}$ at $p\in\Omega$ is the dimension of the orbit of $\mathfrak{N}$ through that point.
We say $\mathfrak{N}$ has positive rank if it has positive rank at every point.

An N-structure is called {\it pure} if the dimension of its stalks is locally constant.
An N-structure is called {\it polarized} if it has positive rank and no singular orbits---this does not mean the orbit dimension (the rank) is locally constant, as the stalk and orbit dimensions may vary together.
An N-structure is called {\it polarizable} if it contains a polarized substructure.
An example of Cheeger-Gromov \cite{CG1} shows the existence of a non-polarizable F-structure on a 4-dimensional manifold.

Let $\Omega$ be a domain that is saturated for some polarized N-structure $\mathfrak{N}$ of positive rank.
An {\it atlas} for $\mathfrak{N}$, denoted $\mathcal{A}=\{(\Omega_i,\,\mathfrak{N}_i)\}_i$, consists of a collection of countably many open sets $\Omega_i$ with $\Omega=\bigcup_i{\Omega}_i$, so that each $\Omega_i$ is saturated under $\mathfrak{N}$ (not just $\mathfrak{N}_i$), so that $\mathfrak{N}_i=(\Omega_i,\mathcal{N}_i,\iota_i)$ is a pure substructure of $\mathfrak{N}|_{\Omega_i}$, and so that the $\Omega_i$ themselves have universal covers $\pi_i:\overline{\Omega}_i\rightarrow{\Omega}_i$ on which the lifted structure $\overline{\mathfrak{N}}_i$ is a constant sheaf whose action integrates to a global action of a connected, simply connected Lie group.
Further, that each $p\in\Omega$ lies in finitely many of the $\Omega_i$, that the stalks $\mathcal{N}_{i,p}$ at $p$ can be ordered by strict inclusion: $\mathcal{N}_{i_1,p}\subset\dots\subset\mathcal{N}_{i_k,p}$, and that there is always some $i$ so that $\mathcal{N}_{i,p}=\mathcal{N}_p$.
Lemma 1.2 of Cheeger-Gromov \cite{CG1} states that an atlas always exists.
An atlas is called {\it polarized} if each pure N-structure $\mathfrak{N}_i$ also has constant rank.

We present some definitions that describe interactions between N-structures and geometry.
A metric is called {\it invariant} under an N-structure if the action of $\mathcal{N}$ is isometric---more precisely, if the image of the monomorphism $\iota:\mathcal{N}\rightarrow\mathcal{X}$ is in the sub-sheaf of Killling fields.
A polarized atlas $\mathfrak{A}=\{(\Omega_i,\,\mathfrak{N}_i)\}$ will be called {\it $C$-regular} if the norm of the second fundamental form of any orbit of $\mathfrak{N}_i$ is bounded from above by $C$, and the multiplicity of the covering $\{\Omega_i\}$ is also bounded by $C$.
A polarized atlas will be called {\it $C$-regular with locally bounded curvature} if the norm of the second fundamental form of any orbit of $\mathcal{N}_i$ at a point $p_i\in{U}_i$ is bounded from above by $Cr_{\RRR}(p_i)^{-1}$, and the multiplicity of the covering $\{U_i\}$ is bounded by $C$.

\begin{lemma}[Global integrability for N-structures] \label{LemmaUniqueGlobalExtensions}
If $\mathfrak{N}$ is any N-structure on a domain $\Omega$ with an invariant metric, and if $\Omega$ is simply connected, then any element $\mathfrak{b}\in\mathcal{N}_p$ of the stalk at any point $p\in\Omega$ extends uniquely to a Killing field $V$ on $\Omega$.
\end{lemma}
{\it Proof}.
On any (differentiable) Riemannian manifold, a Killing field is locally determined by its germ at a point.
In this proof, we shall ``cheat'' slightly by assuming the metric and underlying manifold are analytic (which is true for the manifolds we are interested in).

Suppose $\gamma(t)$, $t\in[0,1]$ is any loop with $\gamma(0)=\gamma(1)=p$.
Assume $\mu_s(t)$ is a homotopy of $\gamma(t)$ to the constant path at $p$; that is, $\mu_0(t)=\gamma(t)$, $\mu_1(t)=p$, and $\mu_s(1)=\mu_s(0)=p$.
We assume that the map $\mu:[0,1]\times[0,1]\rightarrow{N}$ is $C^0$, and is analytic when restricted to $(0,1)\times(0,1)$.

Letting $\mathfrak{v}\in\mathcal{N}_p$, in any sufficiently small neighborhood of $p$, the element $\mathfrak{v}$ has a unique extension to a Killing field $V$.
Covering $\gamma(t)$ by sufficiently small open sets, we obtain an extension of $\mathfrak{v}$ to a Killing field $V$ in some neighborhood of $\gamma(t)$.
This extension is unique along $\gamma$ itself; the issue is that $V(\gamma(0))$ might not equal $V(\gamma(1))$.

Repeating this for any $s\in[0,1]$, we can extend $V$ along the path $t\mapsto\mu_s(t)$.
We obtain a function
\be
\begin{aligned}
\left[0,\,1\right]&\;\longrightarrow\;T_pN \\
s&\;\mapsto\;V(\mu_s(1))
\end{aligned}
\ee
given by $V(\mu_s(1))$.
Because Killing fields on analytic Riemannian manifolds analytic (for instance they satisfy the system $\triangle{V}+\Ric(V)=0$), this map is also analytic.
However, when $s$ is small enough that the path $t\mapsto\mu_s(t)$ lies in a neighborhood of $p$ so small that $V$ is uniquely determined, the map $s\mapsto{V}(\mu_s(1))$ is constant.
By analyticity, it is constant for all $s$.

This shows that given $\mathfrak{v}\in\mathcal{N}_p$, we can define the field $V$ at any point $q$ by connecting $p$ to $q$ with any path and extending $V$ along this path.
The argument above is then used to show that the vector $V(q)$ is independent of the path chosen.
\qed

\begin{theorem}[Cheeger-Gromov \cite{CG2}, Cheeger-Fukaya-Gromov \cite{CFG}] \label{ThmChGrNStruct}
There exists $\tau=\tau(n,\,\delta,\,\alpha)>0$ so that if $\Omega\subset{N}^n$ is a domain in a complete Riemannian manifold $N$ with $|\Riem|<1$ on $\Omega^{(1)}$, and if $Vol\,B(p,1)\le\tau$ for all $p\in\Omega$, then a neighborhood of $\Omega$ exists (that is within $\Omega^{(1)}$) that is saturated with respect to an N-structure $\mathfrak{N}$, and so that the metric on $N$ is $\delta$-close in the $C^{1,\alpha}$ sense to a metric for which $\mathcal{N}$ is invariant.
\end{theorem}
In the case the metric has an elliptic system, $C^{1,\alpha}$-closeness can be improved to $C^{k,\alpha}$-closeness, but where $\tau$ depends also on $k$. 
\begin{theorem}[Cheeger-Rong \cite{CR}] \label{ThmChRong}
If, in addition, $\tau$ is sufficiently small compared to the diameter of $\Omega\subseteq{N}^n$, then $\mathfrak{N}$ is pure.
\end{theorem}
\begin{theorem}[Rong \cite{Rong}] \label{ThmRong}
If, in addition, $\tau$ is sufficiently small and $\Omega\subset{N}^4$ is 4-dimensional, then $\mathfrak{N}$ is polarizable.
\end{theorem}
In addition, there exists a $C<\infty$ so that $\mathfrak{N}$ has a polarized $C$-regular atlas.
\begin{theorem}[Naber-Tian \cite{NT1}] \label{ThmNaberTian}
If $\pi:\Omega\rightarrow\overline{\Omega}$ (where $\Omega\subseteq{N}^4$) is the projection onto the orbit space of a pure N-structure $\mathfrak{N}$, then $\mathfrak{N}$ is an orbifold with $C^\infty$ orbifold points.
\end{theorem}
\begin{theorem}[Cheeger-Fukaya-Gromov \cite{CFG}, Naber-Tian \cite{NT1}, Cheeger-Rong \cite{CR}] \label{ThmCollapsing}
Under the hypotheses of Theorems \ref{ThmChGrNStruct} and \ref{ThmChRong}, and if the metric on $N$ is $\mathfrak{N}$-invariant, the quotient $N\rightarrow\overline{N}$ along the orbits of $\mathfrak{N}$ is a Riemannian orbifold with $C^\infty$ orbifold points, and injectivity radius bounded from below on compact sub-domains.
\end{theorem}

\section{Proof of Theorem \ref{ThmHarnackI}}

\subsection{Outline}

The assertion is that, when energy is small enough, the local scale at $p$ controls the local scale at surrounding points.
If $r_{\mathcal{R}}$ degenerates at some nearby point $p'$, we first re-choose the point $p'$ so that $r_{\mathcal{R}}(p')$ is ``almost'' smallest among all sufficiently nearby $p'$.
Rescaling, we have $r_{\mathcal{R}}(p')=1$ and $r_{\mathcal{R}}$ is bounded uniformly from below on a large region $\Omega$.
With $\int|\Riem|^2$ small, Theorem \ref{ThmCollapsing} forces the existence of an N-structure on $\Omega$.

By passing to the universal cover, we would like the collapsing directions to ``unwrap,'' and become unbounded.
But this is not immediately clear: the manifold could resemble a 3-sphere crossed with a line, where collapse is along Hopf fibers; this is simply connected so passage to the universal cover changes nothing.
But in our situation a Cheeger-Gromoll style splitting theorem implies that the limit is indeed one-ended, so we rule out behavior like $\mathbb{S}^3\times\mathbb{R}$.
Specifically, we prove that collapsing directions must carry homology.

Passing to the universal cover $\widetilde\Omega\rightarrow\Omega$, we know that all orbits of the N-structure are unbounded, and the injectivity radius is bounded from below.
Further, Lemma \ref{LemmaUniqueGlobalExtensions} implies the N-structure is represented by universally defined Killing fields.
The domains $\widetilde\Omega$ then converge to a complete Ricci-flat manifold with $r_{\mathcal{R}}=1$ somewhere, and with at least one Killing field.
Lemma \ref{LemmaFlatMfld} implies that these are flat, contradicting that $r_{\mathcal{R}}$ is not infinite.

\subsection{Point reselection, and properties of the sequence of counterexamples}

In effect, the quantity $\Lambda$ is the scale, and the metric can be rechosen to make $\Lambda=1$; however we will leave $\Lambda$ explicit.
Having chosen $K$, assume there is a sequence of counterexamples, so choose a sequence $\epsilon_i\searrow0$ and for each $\epsilon_i$, a small $\delta_i>0$ and a pointed manifold $(M_i^4,p_i)$ exist so that
\begin{itemize}
\item On $M_i^4$ we have $\Ric\ge-\Lambda^2$.
If the metric is extremal-K\"ahler then $|\nabla{R}|<\Lambda^3$.
If the metric is not CSC, but has another elliptic system, then $|\Ric|<\Lambda^2$.
\item $\int_{B(p_i,\,2Kr_{\mathcal{R}}(p_i))}|\Riem|^2<\epsilon_i$
\item Points $p'_i\in{B}(p_i,\,Kr_{\mathcal{R}}(p_i))$ exist with $r_{\mathcal{R}}(p'_i)/r_{\mathcal{R}}(p_i)<\delta_i$
\item $\delta_i\searrow0$.
\end{itemize}

The first step is to improve the choice of $p_i$ in order to make $r_{\mathcal{R}}$ ``almost'' minimal among all nearby points.
For the moment, we drop $i$ from the notation, so we have a point $p'\in{B}(p,Kr_{\mathcal{R}}(p))$ with $r_{\mathcal{R}}(p')<\delta\,{r}_{\mathcal{R}}(p)$ and $\int_{B(p,Kr_{\mathcal{R}}(p))}|\Riem|^2<\epsilon$.
To initiate the point-picking procedure, let $p_1\in{B}(p',K\delta^{-1}r_{\mathcal{R}}(p'))$ be any point so $r_{\mathcal{R}}(p')<\frac12r_{\mathcal{R}}(p')$, if any such a point exists.
Continuing inductively, having selected points $p'=p_0,p_1,\dots,p_n$, next select a point $p_{n+1}$ in $B(p_n,\,K\delta^{-1}r_{\mathcal{R}}(p_n))$ with $r_{\mathcal{R}}(p_{n+1})<\frac12{r}_{\mathcal{R}}(p_n)$ if such a point exists.
This process terminates with a finite sequence of points $\{p_0,\dots,p_N\}$ so that the final point $p_N$ has the property that $r_{\mathcal{R}}>\frac12r_{\mathcal{R}}(p_N)$ on $B(p_N,\,K\delta^{-1}r_{\mathcal{R}}(p_N))$.

We wish to show that $\int_{B(p_N,\delta^{-1}Kr_{\mathcal{R}}(p_N))}|\Riem|^2<\epsilon$.
Using $r_{\mathcal{R}}(p_j)\le2^{-j}r_{\mathcal{R}}(p_0)$, if $x\in{B}(p_N,\delta^{-1}Kr_{\mathcal{R}}(p_N))$ we have
\be
\begin{aligned}
\dist(x,\,p)
&\;<\;\dist(x,\,p_N)\,+\,\dist(p_N,p_{N-1})\,+\,\dots\,+\,\dist(p_1,p_0)\,+\,\dist(p_0,p) \\
&\;\le\;\delta^{-1}Kr_{\mathcal{R}}(p_{N})\,+\,\delta^{-1}Kr_{\mathcal{R}}(p_{N-1})\,+\,\dots\,+\,\delta^{-1}Kr_{\mathcal{R}}(p_{0})\,+\,Kr_{\mathcal{R}}(p) \\
&\;\le\;\delta^{-1}K\sum_{j=0}^{N}2^{-j}r_{\mathcal{R}}(p_0) \,+\,Kr_{\mathcal{R}}(p).
\end{aligned}
\ee
Since $r_{\mathcal{R}}(p_0)<\delta{r}_{\mathcal{R}}(p)$, we have $\dist(x,\,p)\;<\;2Kr_{\mathcal{R}}(p)$.
Therefore $B\left(p_N,\,\delta^{-1}Kr_{\mathcal{R}}(p_N)\right)\;\subset\;B(p,\,2Kr_{\mathcal{R}}(p))$, so indeed $\int_{B(p_N,\delta^{-1}Kr_{\mathcal{R}}(p_N))}|\Riem|^2<\epsilon$.
Now let the new $p$ be this $p_N$ that we have found.

Re-introduce $i$ into the notation, and re-scale so that $r_{\mathcal{R}}(p_i)=1$.
Labelling the original point $P_i$, the newly chosen $p_i$ has $r_{\mathcal{R}}(p_i)<\delta_i{r}_{\mathcal{R}}(p_i)$.
In the constant scalar curvature case, we have $|R|\le12{r}_{\mathcal{R}}(p_i)$, so in the rescaled metric we have $|R|<12\delta_i^{-1}$.
In the extremal K\"ahler case, where $|\nabla{R}|<\Lambda^3$, we have that, in the rescaled matric, $|R|<C(\delta_i,K)\Lambda^2$ where $\lim_{\delta_i\rightarrow0}C(\delta_i,K)=0$.
The $p_i$ satisfy
\begin{itemize}
\item[{\it{i}})] $\int_{B\left(p_i,\,\frac12K\delta^{-1}_i\right)}|\Riem|^2<\epsilon_i$
\item[{\it{ii}})] On $B\left(p_i,\,\delta^{-1}_iK\right)$, $\Ric\ge-\delta_i^2\Lambda^2$ and $|R|<C(\delta_i,K)\Lambda^2$ where $\lim_iC(\delta_i,K)=0$.
\item[{\it{iii}})] $r_{\mathcal{R}}(p_i)=1$, so a nearby point has $|\Riem|=1$.
\item[{\it{iv}})]  $r_{\mathcal{R}}\;\ge\;\frac12$ so $|\Riem|\le4$ on $B\left(p_i,\,\frac12K\delta_i^{-1}\right)$
\item[{\it{v}})] $\delta_i\searrow0$ and $\epsilon_i\searrow0$.
\end{itemize}
By ({\it{i}}) through ({\it{iv}}) and Lemma \ref{LemmaChTiForcedCollpasing}, the sequence of manifolds-with-boundary $\{B\left(p_i,\,\frac12K\delta_i^{-1}\right)\}_i$ collapses with bounded curvature.
By Theorem \ref{ThmCollapsing} and the comment immediately after, we have a polarized, $C$-regular N-structure $\mathfrak{N}_i$ on a saturation of, say, $B(p_i,\,\frac12K\delta_i^{-1}-1)$.
Let the manifold-with-boundary $N_i\subset{B}(p_i,\,\frac12K\delta_i^{-1})$ be this saturation.

\subsection{Construction of a limiting object, and regularity}

\underline{\it Step I}: {\it Construction and properties of the Gromov-Hausdorff limit $N_\infty$ of the $N_i$}.

By point ({\it{ii}}), Ricci curvature is almost non-negative and scalar curvature converges to zero.
Without any theory of collapse with bounded curvature, Gromov's pre-compactness theorem implies the pointed manifolds $(N_i,p_i)$ converge to a limiting pointed length space $(N_\infty,p_\infty)$ (after passing to a subsequence).
By the Naber-Tian result, Theorem \ref{ThmNaberTian}, $N_\infty$ is a $C^\infty$ orbifold with a Riemannian metric and bounded sectional curvature.

The limit should have $\Ric\ge0$ in a generalized sense, and because $N_\infty$ has dimension 3 or less, this would imply a generalized non-negativity of the curvature operator.
One may expect that the limiting manifold is scalar-flat (so Ricci-flat) in some generalize sense, but the way this passes to the limit is not particularly clear, and we have to work a little harder to show the limits are flat.
The first step is to show limits are one-ended.

\noindent\underline{\it Step II}: {\it The length space $N_\infty$ is one-ended}.

If not, then there exists a line $\gamma_\infty$ on the length space $N_\infty$.
From Cheeger-Colding theory (Theorem 6.64 of \cite{CCo}) the limit $N_\infty$ must have the metric structure of $\mathbb{R}\times{X}_\infty$ for some length space $X_\infty$.
In our situation, with a sectional curvature bound, we arrive at the stronger conclusion that $X_\infty$ is a flat manifold, and that sectional curvature converges pointwise to 0.\footnote{Cheeger-Colding theory is not really necessary, but provides the function $b_\infty$ and shortens the argument somewhat.}

First we re-select the sequence $N_i$.
Choose an exhaustion $\Omega_i$ of $N_\infty$, so that $\Omega_i$ is a domain that satisfies the following three criteria: $\Omega_i$ is connected and has at least 2 ends, $\Omega_i$ contains $B(p_\infty,\,2^i)$, and $\Omega_i$ contains at least two boundary components that are separated by a large distance, say $2^i$.
Now let the new $N_i$ be subsets of the old $N_i$ that are saturated and $2^{-i}$-close to $\Omega_i$ in the Gromov-Hausdorff sense.
We have the projections $\pi_i:{N}_i\rightarrow{N}_\infty$ that collapse the orbits of $\mathfrak{N}_i$ to points.
The map $\pi_i$ is both differentiable and a $2^{-i}$ Gromov-Hausdorff approximation, and the metric $g_i$ on $N_i$ is $2^{-i}$-close in the $C^{k,\alpha}$ sense to a metric for which $\pi_i$ is generically a Riemannian submersion.

Denote by $p_i$ the basepoint in $N_i$, where $p_i\rightarrow{p}_\infty$ as $N_i\to{N}_\infty$.
Since $N_i$ has (at least) two boundary components separated by at least $2^i$, there is a unit-parametrized geodesic path $\gamma_i$ of length at least $2^i$ between them, that lies a uniformly finite distance from $p_i$.

Let $t_i^\pm$ be values so $\gamma_i(t_i^{\pm})$ are the endpoints of $\gamma_i$, and assume that, under the Gromov-Hausdorff approximations $\pi_i:N_i\rightarrow{N}_\infty$, the endpoints $\gamma_i(t_i^{\pm})$ map to the endpoints $\gamma_\infty(t^{\pm}_\infty)$.

Let $b_\infty:N_\infty=\mathbb{R}\times{X}_\infty\rightarrow\mathbb{R}$ be the projection; this is a buseman function associated to the line $\gamma_\infty$.
We may assume $b_\infty(p_\infty)=0$.
Abusing notation, we will also use $b_\infty$ to indicate the pullback functions $\pi_i^{*}(b_\infty)$ on $N_i$.
Since $\pi_i$ is a smooth, almost-Riemannian submersion, $b_\infty$ has uniform $C^{k,\alpha}$ control on $N_i$, based on the characteristics of the submersion.
In particular, the gradient is pinched: $||\nabla{b}_\infty|_{g_i}-1|\le{2}^{-i}$.
We have Hessian bounds but not Hessiam pinching, as it a priori depends on the second fundamental forms of the submersion fibers.

We will use the pointwise pinching of the Ricci tensor and some standard theory to obtain the pinching.
Associated to $\gamma_i$, we have the usual almost-buseman functions:
\be
b_i^{\pm}(x)\;=\;t_i^{\pm}\,-\,\dist(x,\,\gamma_i(t_i^{\pm})).
\ee
We have that that $b_\infty(x)$ is $2^{-i}$-close to $b_i^+$ and to $-b_i^-$ (so also $|b_i^++b_i^-|\le2^{1-i}$).
By the usual Laplacian  comparison argument, we have
\be
\triangle{b}_i^{\pm}\;\ge\;-{3}\cdot{2^{2i}} \label{EstAlmostBuseman}
\ee
in the barrier sense.
Now lift to a local cover $\widetilde{U_i}\rightarrow{U}_i$ so that $p_i\in{U}_i$, where the action of the N-structure $\mathfrak{N}_i$ is generated by connected Lie groups, and the injectivity radius on $\widetilde{U}_i$ is bounded from below.
Lifting $b_\infty$, $b_i^+$, $b_i^-$ to $\widetilde{U_i}$, we retain (\ref{EstAlmostBuseman}), as well as the pointwise estimate $|b_\infty\mp{b}_i^{\pm}|<2^{-i}$

Taking the limit as $i\rightarrow\infty$, we have convergence $\widetilde{U_i}\rightarrow\widetilde{U}_\infty$ in the $C^{k,\alpha}$-sense, where $\widetilde{U}_\infty$ is a 4-dimensional manifold-with-boundary.
Then $b_i^\pm$ converge to functions $b_\infty^\pm$ on $\widetilde{U}_\infty$ where $b_\infty^+=-b_\infty^-=b_\infty$, and where the $\triangle{b}_\infty^\pm\ge0$ in the barrier sense.

Thus $b_\infty$ on $\widetilde{U}_\infty$ is harmonic, and since $\widetilde{U}_\infty$ has $\Ric=0$, the B\"ochner formula gives $|\nabla^2b_\infty|^2=0$.
Of course $b_\infty$ is a distance function, so the domain $\widetilde{U}_\infty$ has a metric splitting into a Ricci-flat 3-manifold and a line segment; thus it is flat.
This contradicts point ({\it{iii}}) above, which says that $r_{\mathcal{R}}(p_i)=1$, a condition that lifts to $\widetilde{U_i}$ and that passes to the limit $\widetilde{U}_\infty$.
Thus indeed $N_\infty$ to be one-ended.

\noindent\underline{\it Step III}: {\it Reselection of the domains $N_i$, and proof that $\chi(N_i)\ge1$}.

As in the previous step, let $\Omega_i\subset{N}_\infty$ be an exhaustion of $N_\infty$ by connected, pre-compact domains.
Let the new $N_i$ be a saturated subset of the old $N_i$ that is $2^{-i}$ close in the Gromov-Hausdorff sense to $\Omega_i$.
We have smooth projections $\pi_i:N_i\rightarrow\Omega_i$ that collapse orbits of the N-structures  $\mathfrak{N}_i$ to points, and we may assume the metrics $g_i$ on $N_i$ are $2^{-i}$-close in the $C^{k,\alpha}$ sense to invariant metrics, so $\pi_i$ is generically an almost-Riemannian submersion.
Finally, pass to the universal cover, so $N_i$ is simply connected.
We have
\begin{itemize}
\item The basepoint $p_i\in{N}_i$ has $r_{\mathcal{R}}(p_i)=1$
\item $\Ric\ge-\delta_i^2\Lambda^2$ and $|R|\le{C}(\delta_i,K)\Lambda^2$ on $N_i$, where $\lim_iC(\delta_i,K)=0$
\item $N_i$ is saturated by a pure, $C$-regular N-structure $\mathfrak{N}_i$ of positive rank
\item $N_i$ is connected and simply connected
\item $N_i$ is one-ended
\end{itemize}

\begin{lemma} \label{LemmaEulerNumbers}
The Euler number $\chi$ on $N_i$ is strictly positive.
\end{lemma}
{\it Proof}.
By simple connectedness, we have $H^1(N_i;\mathbb{R})=\{0\}$, so by Poincare duality $H^3(N_i,\partial{N}_i;\mathbb{R})=\{0\}$.
Then the relative homology sequence gives
\be
\dots
\rightarrow{H}^3(N_i,\partial{N}_i;\mathbb{R})
\rightarrow{H}^3(N_i;\mathbb{R})
\rightarrow{H}^3(\partial{N}_i;\mathbb{R})
\rightarrow{H}^4(N_i,\partial{N}_i;\mathbb{R})
\rightarrow{H}^4(N_i;\mathbb{R})
\rightarrow\dots
\ee
Now $H^4(N_i,\partial{N}_i;\mathbb{R})=\mathbb{R}$ (generated by the fundamental class) and $H^4(N_i;\mathbb{R})=\{0\}$.
By one-endedness $H^3(\partial{N}_i;\mathbb{R})=\mathbb{R}$, so exactness forces $H^3(N_i;\mathbb{R})=\{0\}$.
Therefore the Euler number of $N_i$ is $\chi(N_i)=1+b^2(N_i)\ge1$ (where $b^2(N_i)$ is the second betti number of $N_i$).
\qed

By Lemma \ref{LemmaUniqueGlobalExtensions}, any element of a stalk $\mathfrak{v}\in\mathcal{N}_p$ has a unique global extension to a vector field $V$.
Because the metric on $N_i$ is close to a metric for which $V$ is Killing, the integral curves of $V$ are either all unbounded, or all bounded.
If the integral curves of such a $V$ were both bounded and no-where zero, it would force $\chi(N_i)=0$, contradicting this lemma.

\subsection{The contradiction argument}

Passing to a subsequence, we can assume the stalks of the sheafs $\mathcal{N}_i$ on the manifold $N_i$ have constant dimension.
We complete the argument by ruling out the three possible dimensions.

\begin{lemma} \label{LemmaOneDimStructure}
The sheaves $\mathcal{N}_i$ do not have stalk dimension $1$.
\end{lemma}
{\it Proof}.
If $\mathfrak{N}_i$ had rank 1, it would be represented by a global vector field $V$, which has no zeros (by positivity of rank).
It is impossible that its orbits are bounded or else $\chi(N_i)=0$; thus its orbits are unbounded.
Because the $N_i$ have no further collapsing directions, the injectivity radii are bounded from below.
Taking a limit as $N_i\rightarrow{N}_\infty$, we have that $N_\infty$ is a complete Ricci-flat 4-manifold with a Killing field $V$ that has a uniform bound on $|\nabla{V}||V|^{-1}$ (by $C$-regularity), and has $|\Riem|=1$ at at least one point.
By Lemma \ref{LemmaFlatMfld}, such a manifold is flat, a contradiction.
\qed

\begin{lemma}
The $\mathcal{N}_i$ do not have stalk dimension $2$.
\end{lemma}
{\it Proof}. If $\mathfrak{N}_i$ has rank 2, the sheaf $\mathcal{N}_i$ is abelian.
We may choose two vector fields $V_1$ and $V_2$ to represent $\mathcal{N}_i$, and we may assume the closures of the orbits of either field is 1-dimensional.
Further, by $C$-regularity of the N-structure, we may assume that either $|\nabla{V}_1||V_1|^{-1}$ is uniformly bounded, or else $V_1$ has a zero, and likewise for $V_2$.

If either of the $V_j$ ($j=1,2$) has orbits that are bounded, then $V_j$ must have a zero somewhere, or else it would force the absurdity $\chi(N_i)=0$.

Assume $V_1$, say, has a zero and $V_2$ does not.
Then the orbits of $V_2$ are unbounded and $|\nabla{V}_2||V_2|^{-1}$ is bounded.
Near the zero-set of $V_1$, the injectivity radius is bounded from below, or else, if collapse still occurred, there would be a third dimension to the stalks of $\mathfrak{N}_i$.
Now this zero might be increasingly far away as $i\rightarrow\infty$, in which case we can shrink the domains $N_i$ somewhat, and we are in the previous case, so we can assume the zero-set is a finite distance from $p_i$.
Since the injectivity radius is bounded a finite distance from $p_i$, it is bounded on compact sets containing $p_i$, and again we get convergence to a complete Ricci-flat 4-manifold $N_\infty$, that has $|\Riem|=1$ at at least one point, and has a nowhere-zero Killing field $V_2$, again contradicting Lemma \ref{LemmaFlatMfld}.

Finally, if $|V_1|$ and $|V_2|$ both have zeros on $N_i$, then all orbits are bounded.
The zero-sets are non-intersecting by positivity of rank.
Thus $\chi(N_i)=0$ by an easy argument (eg. in Proposition 1.5 of \cite{CG1}), and again we have a contradiction.
\qed

\begin{lemma}
The $\mathcal{N}_i$ do not have stalk dimension $3$.
\end{lemma}
{\it Proof}.
Quotienting out by the orbits of $\mathfrak{N}_i$, we have a pointed limit $(\pi_i(N_i),\pi_i(p_i))\rightarrow(N_\infty,p_\infty)$, where $|\Riem|=1$ at $p_i$ and $|\Riem|<4$ on $N_i$.
Because $\mathfrak{N}_i$ as rank 3, $N_\infty'$ is either a line or a ray.
The first case is 2-ended, so is impossible.
In the case of a ray, there is a single singular orbit, $\mathcal{O}_i$ in each $N_i$, over the ray point.
This must be a singular orbit, for if it is non-singular it is 3-dimensional and therefore separates $N_i$, so its image under $\pi_i:N_i\rightarrow{N}_\infty$ is an orbifold point which elsewhere in the ray, which is impossible.

Note also that no exceptional orbits exist, except possibly on the singular orbit itself; this is again seen by noting that $k$-to-$1$ exceptional orbits translate to $k$-to-$1$ orbifold points on the quotient, which has no orbifold points except the ray point.
Thus with neither singular nor exceptional orbits except $\mathcal{O}_i$, the manifold $N_i$ has a deformation ratract onto $\mathcal{O}_i$.
Thus the submanifold $\mathcal{O}_i$ is either a simply connected 1-manifold with an almost-Killing field, and therefore $\mathbb{R}^1$, or a simply-connected 2-manifold with two commuting, no-where zero, almost-Killing fields, and is therefore $\mathbb{R}^2$ with and almost-flat metric.

Because the singular orbit $\mathcal{O}_i\subset{N}_i$ is non-collapsed, the injectivity radius of $N_i$ near $\mathcal{O}_i$ is bounded away from zero.
By the boundedness of curvature, the injectivity radius of $N_i$ is bounded away from zero on compact domains, so we get convergence $N_i\rightarrow{N}_\infty$ to a complete Ricci-flat 4-manifold $N_\infty$ with $|\Riem|=1$ at at least one point.
Further $N_\infty$ has three Killing fields, at least one of which has unbounded orbits.
By Lemma \ref{LemmaFlatMfld} again, $N_\infty$ is flat, a contradiction.
\qed

Having shown that the pointed manifolds $(N_i,p_i)$ are collapsed with bounded curvature, but that the resulting N-structure does not have rank 1, 2, or 3 (or obviously 4), we conclude that it is impossible to find such manifolds $N_i$, contradicting that collapse happens.
This establishes the theorem.

\section{Proof of Theorem \ref{ThmHarnackII} and Corollary \ref{CorEllipticEsts}}

We restate Theorem \ref{ThmHarnackII} for convenience: \\

{\it 
Given $K<\infty$, $k>0$, and $\mu\in(0,1]$, there exist numbers $\epsilon_0=\epsilon_0(K)>0$ and $C=C(\mu,k,K)<\infty$ so that the following holds.
If $r_{\mathcal{R}}(p)<K\Lambda^{-1}$ and  $\int_{B(p,\,2{r}_{\mathcal{R}}(p))}|\Riem|^2\;\le\;\epsilon_0$, then
\be
r_{\mathcal{R}}(p)^{-2k}\;\le\;\frac{C}{\Vol\,B(p,\,\mu{r}_{\mathcal{R}}(p))}\int_{B(p,\,\mu{r}_{\mathcal{R}}(p))}|\Riem|^{k}. \label{IneqReverseHarnackII}
\ee
}

Fix $\Lambda,\mu,k$, and assume there is no such $C$, meaning there is a sequence of counterexamples so that the quantities
\be
\frac{r_{\mathcal{R}}(p_i)^{2k}}{\Vol\,B(p_i,\,\mu{r}_{\mathcal{R}}(p_i))}\int_{B(p_i,\,\mu{r}_{\mathcal{R}}(p_i))}|\Riem|^{k}
\ee
can degenerate to zero no matter what $\epsilon_0>0$ is chosen.
By Theorem \ref{ThmHarnackI}, we can choose $\epsilon_0$ small enough that there is a $\delta_0$ so that $r_{\mathcal{R}}\ge\delta_0r_{\mathcal{R}}(p_i)$ on $B(p_i,\,2r_{\mathcal{R}}(p_i))$.
Then the exponential map has no conjugate points on some ball of radius definitely (though slightly) larger that $r_{\mathcal{R}}(p_i)$.
Namely $\exp_{p_i}:B(o_i,(1+\eta)r_{\mathcal{R}}(p_i))\rightarrow{B}(p_i,\,(1+\eta)r_{\mathcal{R}}(p_i))$ is a local homeomorphism, where $\eta$ is independent of $i$, and $o_i$ is the origin in the tangent space at $p_i$.
Lifting to the tangent space at $p_i$, we have a ball $B((1+\eta)r_{\mathcal{R}}(p_i))$ that is contractible.
Finally scale so that $r_{\mathcal{R}}(p_i)=1$.

Now the exponential map $exp_{p_i}:B(o_i,1+\eta)\rightarrow{B}(p_i,1+\eta)$ does not evenly cover the target.
To get around this, choose an open fundamental domain in the usual way: the unique pre-image under $exp_{p_i}$ of $B(p_i,1+\eta)$ minus the cut locus, that contains the origin.
By choosing different pre-images of $p_i$, the fundamental domain withing the tangenst space is also shifted.
If the chosen pre-image of $p_i$ is within a certain definite distance, say $d_0$, of the origin, the exponential map will still have no conjugate points, so the exponential map will be a local diffeomorphism, and will be 1-1.

Let $\Omega_i$ be the interior of the closure of the union of ``fundamental domains'' in the tangent space at $p_i$ whose basepoints are a distance less than $p_0$ from the origin.
Then $exp_{p_i}$ restricted to $\Omega_i$ is an even covering; it is, say, a $M$-to-$1$ cover where $M\in\mathbb{N}$ depends on $i$.
To account for the $\mu$, let $\Omega^\mu_i\subset\Omega_i$ be the union of the subsets of the same fundamental domains, restricted to the various pre-images of $B(p_i,\mu)$ instead of all of $B(p_i,1+\eta)$.
Giving $\Omega_i$ the pullback metric, we have
\be
\frac{1}{\Vol\,\Omega^\mu_i}\int_{\Omega^\mu_i}|\Riem|^k
\;=\;
\frac{1}{M\Vol\,B(p_i,\mu)}\cdot{M}\int_{B(p_i,\mu)}|\Riem|^k \;\stackrel{i\rightarrow0}{\longrightarrow}\;0.
\ee
Now we have that $|\Riem|=1$ somewhere on $B(p_i,1)\subset{B}(p_i,1+\eta)$ and therefore on $\Omega^1_i\subset\Omega_i$.
Taking a limit as $i\rightarrow\infty$, the Riemannian domains converge $\Omega_i\rightarrow\Omega_\infty$ in the $C^\infty$ sense and the limiting metric on $\Omega_\infty$ is Ricci-flat, has $|\Riem|=1$ somewhere on its interior.
Letting $\Omega^\mu_\infty=\lim_i\Omega^\mu_i$, we also have
\be
\frac{1}{\Vol\,\Omega^\mu_\infty}\int_{\Omega^\mu_\infty}|\Riem|^k\;=\;0.
\ee
By the classic Harnack inequality for elliptic systems, this is impossible.

Corollary \ref{CorEllipticEsts} is proved similarly.
Choose $l$, pick counterexamples $B(p_i,r_{\mathcal{R}}(p_i))$, and scale so $r_{\mathcal{R}}(p_i)=1$.
Again passing to the tangent spaces of the $p_i$, we have convergence of the metrics on the slightly larger, contractible manifolds $B(o_i,1+\eta)$.
The limiting metric on $B(o_\infty,1+\eta)$ has definite bounds on the quantities $|\nabla^l{r}_{\mathcal{R}}|$ within $B(o_\infty,1+\frac12\eta)$, so by $C^\infty$ convergence, the stated bounds must hold on the $B(o_i,1+\frac12\eta)$, and so on the original $B(p_i,r_{\mathcal{R}}(p_i))$.


\begin{thebibliography}{9}

\bibitem{And4} {M. Anderson}, \emph{Canonical metrics on 3-manifolds and 4-manifolds}, {Asian Journal of Mathematics}. Vol. 10, No. 1 (2006) {127--164}

\bibitem{CCo} {J. Cheeger and T. Colding}, \emph{Lower bounds on Ricci curvature and the almost rigidity of warped products}. {Annals of Mathematics} Vol. 144, No. 1 (1996) {189--237}

\bibitem{CFG} {J. Cheeger, K. Fukaya and M. Gromov}, \emph{Nilpotent structures and invariant metrics on collapsed manifolds}. {Journal of the American Mathematical Society}, Vol. 5, No. 2 (1992) {327--372}

\bibitem{CG1} {J. Cheeger and M. Gromov}, \emph{Collapsing Riemannian manifolds while keeping their curvature bounded I}. {Journal of Differential Geometry}, Vol. 23, No. 3 (1986) {309--346}

\bibitem{CG2} {J. Cheeger and M. Gromov}, \emph{Collapsing Riemannian manifolds while keeping their curvature bounded II}. {Journal of Differential Geometry}, Vol. 32, No. 1 (1990) {269--298}

\bibitem{Derd} {A. Derdzi\'nski}, \emph{Riemannian Manifolds with Harmonic Curvature}. {Lecture Notes in Mathematics 1156}, Springer-Verlag, Berlin, Heidelberg, New York, 1985

\bibitem{Rong} {X. Rong}, \emph{The existence of polarized F-structures on volume collapsed 4-manifolds}. {Geometrics and Functional Analysis}, Vol. 3, No. 5 (1993) {474--501}

\bibitem{CR} {J. Cheeger and X. Rong}, \emph{Existence of polarized F-structures on collapsed manifolds with bounded curvature and diameter}. {Geometric and Functional Analysis}, Vol. 6, No. 3 (1996) {411--429}

\bibitem{CT} {J. Cheeger and G. Tian}, \emph{Curvature and injectivity radius estimates for Einstein 4-manifolds}. {Journal of the American Mathematical Society}, Vol. 12, No. 2 (2005) {487-525}

\bibitem{NT1} {A. Naber and G. Tian}, \emph{Geometric structures of collapsing Riemannian manifolds, I}. {arXiv:0804.2275}

\bibitem{Web2} {B. Weber}, \emph{Energy and Asymptotics of Ricci Flat 4-Manifolds with a Killing Field}. 

\end{thebibliography}
\end{document}